\documentclass[12pt]{amsart}
\usepackage{amssymb}

\newtheorem{rmk}{Remark}
\newtheorem{thm}{Theorem}

\newtheorem{lem}{Lemma}

\newcommand{\inner}[2]{\mbox{$\mathcal{h} #1,#2 \mathcal{i}$}}

\newcommand{\BR}{\mbox{$\mathbf{R}$}}

\newcommand{\Lg}{\mbox{$\mathfrak{g}$}}
\newcommand{\Lk}{\mbox{$\mathfrak{k}$}}
\newcommand{\Lp}{\mbox{$\mathfrak{p}$}}
\newcommand{\Ls}{\mbox{$\mathfrak{s}$}}
\newcommand{\La}{\mbox{$\mathfrak{a}$}}
\newcommand{\Lm}{\mbox{$\mathfrak{m}$}}
\newcommand{\Lb}{\mbox{$\mathfrak{b}$}}

\newcommand{\Ad}{\mbox{Ad}}
\newcommand{\ad}{\mbox{ad}}

\newcommand{\Pf}{{\em Proof}. }
\newcommand{\EPf}{\hfill$\square$}

\begin{document}


\title{\bf A CLASS OF COMPLETE EMBEDDED MINIMAL SUBMANIFOLDS\\
 IN NONCOMPACT SYMMETRIC SPACES}


\author{Claudio Gorodski}

\date{\today} 

\address{Instituto de Matem\' atica e Estat\'\i stica\\
         Universidade de S\~ ao Paulo\\
         Rua do Mat\~ao, 1010\\
         S\~ ao Paulo, SP 05508-900\\
         Brazil}
\email{\texttt{gorodski@ime.usp.br}}

\begin{abstract}
We show that a totally geodesic submanifold of a symmetric space
satisfying certain conditions admits an extension to a minimal 
submanifold of dimension one higher, and we apply this result to 
construct new examples of complete embedded minimal submanifolds in
simply connected noncompact globally symmetric spaces. 
\end{abstract}

\maketitle 


\section{Introduction}

A submanifold of a Riemannian manifold 
is called \emph{totally geodesic} if its
second fundamental form identically vanishes, and it is called
\emph{minimal} if the trace of its second fundamental form 
identically vanishes. 
Of course, any totally geodesic submanifold is 
minimal. In this paper, in the presence of enough symmetries
in the ambient space, we give a sufficient 
condition to be able to ``extend'' a given 
totally geodesic submanifold to a minimal submanifold 
of dimension one higher, namely to 
construct a new minimal submanifold which contains the
original totally geodesic submanifold as a hypersurface of codimension one. 

In fact, suppose that $S$ is a totally geodesic submanifold 
of a symmetric space $M$ of codimension greater than one.
Let $\gamma$ be a geodesic in $M$ which is
normal to $S$ at $\gamma(0)\in S$. The transvections along
$\gamma$ define a (local) one-parameter group $\{\psi_t\}$ of 
(local) isometries of~$M$. We describe a condition on 
the pair~$(S,\gamma'(0))$ which guarantees that the submanifold obtained by 
translation of $S$ using~$\{\psi_t\}$ is minimal
(and in general not totally geodesic). Under reasonable 
global conditions on $M$ and $S$, 
this result is then extended to produce a minimal embedding
of~$\BR\times S$ into $M$. 

The construction can be applied to a variety of totally geodesic
submanifolds in symmetric spaces, yielding (mostly new) 
examples of minimal submanifolds.
For example, there is a 
rather large class of totally geodesic
submanifolds in symmetric spaces, the
so-called \emph{reflective submanifolds} (which were introduced 
and classified by Leung in~\cite{Le1,Le2,Le3}), which can be translated 
along \emph{any} normal geodesic to produce a minimal submanifold. 
It is worth here noting the case of complex hyperbolic space
$\mathbf{H}^n(\mathbf{C})$. 
Up to ambient isometries, the maximal 
reflective submanifolds in~$\mathbf{H}^n(\mathbf{C})$ turn out to be
of two types:
complex totally geodesic (isometric to
$\mathbf{H}^{n-1}(\mathbf{C})$) or
totally real totally geodesic (isometric to $\mathbf{H}^n(\mathbf{R})$). 
It is easy to show that the minimal hypersurface 
constructed by 
translating a codimension two complex totally geodesic submanifold
along a normal geodesic is a \emph{bisector}.
Bisectors in~$\mathbf{H}^n(\mathbf{C})$ 
are defined as the geometric loci of points equidistant 
from two distinct given points, they are all congruent and
have been used in the construction of 
fundamental polyhedra for discrete groups of isometries 
of $\mathbf{H}^n(\mathbf{C})$ (see~\cite{Go} for these concepts). 
On the other hand, the 
translation of a totally real totally geodesic submanifold 
along a normal geodesic produces a new congruence class of
minimal submanifolds of $\mathbf{H}^n(\mathbf{C})$.
In the case~$n=2$ these are minimal hypersurfaces, and 
are expected as well to play an important role in the construction of 
discrete groups of isometries of $\mathbf{H}^2(\mathbf C)$.
The author would like to thank N. Gusevskii for the suggestion to 
consider this hypersurface, which originated the whole subject
of this paper. 

\section{The construction}

\begin{thm}
(a) Let $M$ be a locally symmetric space, $S\subset M$ 
a totally geodesic submanifold of codimension greater than one,
$p\in S$, and $X$ an infinitesimal transvection which extends a nonzero vector
$X_p\perp T_pS$ and whose local one-parameter group of isometries 
is denoted by~$\{\psi_t\}$. 
Suppose that for every infinitesimal transvection
$Y$ at~$p$ we have that
\[ (L_XL_Y^{2n+1}X)_p \in T_p S\quad\mbox{for $n=0,1,2,\ldots$} \]
Then there exist $\epsilon>0$ and
a neighborhood $U$ of $p$ 
in $S$ such that the map
\[ f:(-\epsilon,\epsilon)\times U \to M,\quad f(t,x)=\psi_t(x) \]
is a minimal immersion. 

(b) If, in addition to the above assumptions, 
$M$ is nonpositively curved, globally symmetric 
and simply connected\footnote{Simple connectedness is automatic if $M$ is 
negatively curved and globally symmetric (\cite{Wo}, Cor.~8.3.13).}, 
and $S$ is complete, then the map~$f$
extends to a complete minimal embedding~$\BR\times S\to M$. 
\end{thm}

\begin{rmk}
\em 
Let $\Lp$ be the set of all infinitesimal transvections at~$p$
and $\Lk$ be the Lie algebra of the local isotropy group at~$p$.
Then $\Lg=\Lk+\Lp$ is a Lie algebra of vector fields on a neighborhood of 
$p$ and it is the involutive Lie algebra associated to $M$
at~$p$ (see~\cite{Wo}, section~8.1). Recall that every infinitesimal
transvection at $p$ extends some vector in $T_pM$ and this 
defines an identification of $\Lp$ with $T_pM$; in particular, 
$\Lp$ is equipped with an $\ad_{\Lk}$-invariant inner product.
Let $\Ls\subset\Lg$ correspond to $T_pS\subset T_pM$.
Then the property of $S$ being totally geodesic in $M$ corresponds to 
the fact that $\Ls\subset \Lg$ is a Lie triple system, namely 
$[[\Ls,\Ls],\Ls]\subset\Ls$. Now the assumption in the theorem
can be rewritten as
\begin{equation}\label{cond'}
 [X,\ad_Y^{2n+1}X] \in \Ls \quad\mbox{for $Y\in\Ls$ and $n=0,1,2,\ldots$}
\end{equation}
In applications we shall often refer to this equivalent condition 
in order to check that the hypothesis of the theorem is fulfilled. 
Note also that condition~(\ref{cond'}) is in general too weak 
to imply that $f$ is a totally geodesic immersion, for
$\Ls+\BR X$ need not be a Lie triple system of~$\Lg$.  
\end{rmk}

{\em Proof of the theorem.} 
(a) Let $V$ be a normal neighborhood of $p$ in $M$ where the 
Killing field~$X$ is defined. 
We first show that~$X$ is normal to $S$ in $U=S\cap V$. 
In fact let $q\in U$ be arbitrary and let $c(t)$, $0\leq t\leq1$,
be a geodesic segment in $U$ joining 
$p$ to $q$. Take any parallel 
vector field $Y(t)$ along $c(t)$ which is everywhere tangent to 
$S$ and define the function $f(t)=\inner{X}{Y}_{c(t)}$. 
We have that $f$ is a real analytic function, because $M$ is an
analytic Riemannian manifold (see Proposition~5.5 in~\cite{He}). 
Since~$X$ is orthogonal to~$S$ at~$p$, 
we also have that~$f(0)=0$. Consider 
the linear operator $\mathcal J:T_pM\to T_pM$ defined 
by~$\mathcal J(v)=R(c'(0),v)c'(0)$. 
The restriction of $X$ to a vector field along $c$, $X(t)$,
is a Jacobi field, therefore $X''(0)=\mathcal J(X(0))$.
Next we use the facts that $\mathcal J$ is a symmetric
operator with respect to the inner product in $T_pM$, 
$T_p S$ is a $\mathcal J$-invariant subspace (because $S$ 
is totally geodesic in $M$), $c''=Y'=0$, $\nabla R=0$ (because $M$
is locally symmetric) and $X'(0)=0$ (because $X$ is an 
infinitesimal transvection at~$p$) to compute that 
\[ f^{(2n)}(0)=\inner{{\mathcal J}^n(Y(0))}{X(0)}=0\quad\mbox{and} \]
\[ f^{(2n+1)}(0)=\inner{{\mathcal J}^n(Y(0))}{X'(0)}=0\quad
\mbox{for $n=0,1,2,\ldots$} \]
Since the derivatives of all orders of~$f$ vanish at~$0$, 
it follows that $f$ identically vanishes. Therefore
$\inner{X}{Y}_{c(1)}=f(1)=0$. Since $Y(1)$ can be any vector in $T_qS$,
this proves that $X$ is orthogonal to $S$ at $q$. 

We have that, $\{\psi_t\}$
being a local one-parameter group of isometries, there exists 
$\epsilon>0$ such that $f(t,x)$ is well-defined 
and satisfies $f(t,x)=(\psi_t\circ f)(0,x)$ 
for~$(t,x)\in(-\epsilon,\epsilon)\times U$.
Note that, for any $x\in U$, 
the image of the differential of~$f$ at~$(0,x)$ 
equals the subspace~$\BR X_x+T_xS$, where $X_x$ is normal to~$T_xS$,
so $f$ is an immersion at~$(0,x)$ if and only if~$X_x\neq0$. 
It is then clear that we can shrink $U$, if necessary, 
so that $f$ is an immersion on~$(-\epsilon,\epsilon)\times U$. Moreover, 
in order to see that $f$ is a minimal immersion, it is enough to
compute the mean curvature vector of $f$ for $t=0$. Let $q\in U$ and 
$n$ be a normal unit vector field along $f$ near $f(0,q)=q$. 
If $v\in T_qS$, then $\inner{\nabla_vn}{v}_q=0$ because $n$ is 
normal to $S$ and $S$ is totally geodesic. Since~$X_q\perp T_qS$
and~$X_q\neq0$,
we only need to check that 
$\inner{\nabla_\frac{X}{||X||}n}{\frac{X}{||X||}}_q=
-\frac{1}{||X||^2}\inner{n}{\nabla_XX}_q=0$. 

So next we show that $(\nabla_XX)_q\in T_qS$.
Consider the involutive Lie algebra $\Lg=\Lk+\Lp$ 
associated to $M$ at~$p$ as in the remark.
There exists $Y\in T_pS$ such that 
$q=\exp_pY$. Let~$\varphi$ be the local transvection from $p$ to $q$ along
$t\mapsto\exp_ptY$. Then $\varphi(p)=q$, the 
differential~$d\varphi_p:T_pM\to T_qM$ is parallel transport along 
$t\mapsto\exp_ptY$ and 
maps $T_pS$ onto~$T_qS$. Moreover~$(\nabla_XX)_q=d\varphi_p(\nabla_ZZ)_p$, 
where $Z=\varphi_*^{-1}X$. 
So it is enough to verify that~$(\nabla_ZZ)_p\in T_pS$.

Note that the local flow defined by $Z$ is $\{\varphi^{-1}\psi_t\varphi\}$
and, on the Lie group level, 
$\varphi^{-1}\psi_t\varphi=\exp t\Ad_{\varphi^{-1}}X$,
so
\[ Z=\Ad_{\varphi^{-1}}X=e^{-\ad_Y}\cdot X=Z^{\Lk}+Z^{\Lp}, \]
where 
\[ Z^{\Lk}=\sum_{n=0}^{\infty}\frac{1}{(2n)!}\ad_Y^{2n}X\quad\mbox{and}\quad
Z^{\Lp}=-\sum_{m=0}^{\infty}
\frac{1}{(2m+1)!}\ad_Y^{2m+1}X \]
are respectively the $\Lk$ and $\Lp$ components of $Z\in\Lg$ with respect 
to the decomposition $\Lg=\Lk+\Lp$. According to results in ch.~X,
section~2 and ch.~XI, section 3 of~\cite{K-N} we have that 
$(\nabla_ZZ)_p=[Z^{\Lk},Z^{\Lp}]$. Therefore
\[ (\nabla_ZZ)_p=-\sum_{n=0}^{\infty}\sum_{m=0}^{\infty}\;
    \frac{1}{(2n!)(2m+1)!}\,[\ad_Y^{2n}X,\ad_Y^{2m+1}X]. \]
Our claim that $(\nabla_ZZ)_p\in T_pS$ 
is now a consequence of condition~(\ref{cond'}) and the lemma below. 
This completes the proof that $f$ is a minimal immersion. 

\smallskip

(b) Since here $M$ is globally symmetric, we have that 
the Killing field $X$ is globally defined, and 
$\psi_t$ is defined for all $t\in\BR$ and it is a global isometry,
so $f$ is well-defined 
on $\BR\times S$. We claim that $f$ is an immersion everywhere. 
It suffices to check this at a point $(0,q)\in\BR\times S$. 
Since the restriction of 
$X$ to the geodesic connecting $p$ and $q$ is a nontrivial
Jacobi field and $M$ is nonpositively curved,
we get that $X_q$ is not zero.
This implies that~$f$ is an immersion at~$(0,q)$. 

We next show that~$f$ is injective. 
Note that this is equivalent to having
$\psi_t(S)\cap S=\varnothing$ for~$t\neq0$. 
In the following we prove in fact that~$d(S,\psi_t(S))=|t|||X_p||$, 
where~$d$ is the metric distance on~$M$. 
Let $\gamma$ be the geodesic in~$M$ through~$p$ with initial 
velocity~$X_p$.
Since $S$ is totally geodesic in~$M$ and~$M$ is nonpositively curved, 
the function~$\rho:x\mapsto d^2(x,S)$ is convex on~$M$
(see~\cite{B-O}). But $T$ is totally geodesic, so the 
restriction $\rho|_T$ is also convex. Now 
$\gamma$ is a common perpendicular 
from $S$ to $T$, so a first variation argument shows that 
$\gamma(t)\in T$ is a critical point of $\rho|_T$. 
Thus~$\gamma(t)$ must be global minimum point of~$\rho_T$.
Clearly~$\gamma(t)\not\in S$, and there cannot be 
another perpendicular from~$\gamma(t)$ to~$S$ 
besides~$\gamma$ for otherwise we would have a geodesic
triangle with two right angles in~$M$ which is impossible
since~$M$ is nonpositively curved. Hence
$d(S,T)=d(\gamma(0),\gamma(t))=|t|||X_p||$. 

We finally prove that $f$ is a proper map. Let $\{(t_n,x_n)\}$ be a 
sequence in $\BR\times S$ such that $(t_n,x_n)\to\infty$ 
(If $\{a_n\}$ is a sequence in a topological space $A$, 
the notation $a_n\to\infty$ means that 
the sequence eventually leaves any given compact 
subset of $A$.) We need to prove that~$f(t_n,x_n)\to\infty$. 
Since $\BR\times S$ has the product topology, we have that 
either~$t_n\to\infty$ or~$x_n\to\infty$. 
In the former case, as we have seen above,
$d(p,\psi_t(x))\geq|t|||X_p||$ for~$x\in S$, and this implies that
$f(t_n,x_n)\to\infty$. In the latter case,~$f(t_n,x_n)\to\infty$
is an easy consequence of the fact that for every~$q\in S$
the map~$t\mapsto d(p,\psi_t(q))$ has a global minimum for~$t=0$,
which can be seen as follows. Of course we may assume that~$q\neq p$. 
If this map has a minimum for some~$t\neq0$, 
then the points~$p$, $q$ and~$\psi_t(q)$ specify a 
geodesic triangle with two right angles contradicting the fact 
that~$M$ is nonpositively curved. Therefore the minimum 
is attained for~$t=0$. This completes the proof of the theorem. \EPf

\begin{lem}\label{alg}
Let $Y\in\Ls$. If $[X,\ad_Y^{2m+1}X]\in\Ls$ for all $m\geq0$, then 
$[\ad_Y^{2n}X,\ad_Y^{2m+1}X]\in\Ls$ for all~$n,m\geq0$.
\end{lem}

\Pf We prove by induction on $n+m$ that 
\[ \ad_Y[\ad_Y^{2n}X,\ad_Y^{2m}X]\in\Ls\quad\mbox{and}\quad
   [\ad_Y^{2n}X,\ad_Y^{2m+1}X]\in\Ls. \]
In the initial case of $n+m=0$ we must have $n=m=0$ and then the
assertion is clear. 

Let $n,m\geq0$ be such that $n+m>0$. We first show that 
$\ad_Y[\ad_Y^{2n}X,\ad_Y^{2m}X]\in\Ls$. Without loss 
of generality we may assume that $n\geq m$. We proceed 
by induction on~$n-m$. The case~$n-m=0$ is trivial. 
Consider the case~$n-m=1$ and write, via the Jacobi identity,
\[ \ad_Y[\ad_Y^{2n}X,\ad_Y^{2m}X]=\ad_Y^2[\ad_Y^{2n-1}X,\ad_Y^{2m}X]
-\ad_Y[\ad_Y^{2n-1}X,\ad_Y^{2m+1}X]. \]
The first summand of the right hand side belongs to $\Ls$ by the induction 
hypothesis on $n+m$ (and the fact that $\ad_Y^2\Ls\subset\Ls$,
since $\Ls$ is a Lie triple system) and the second summand is zero
because $n-m=1$. Therefore $\ad_Y[\ad_Y^{2n}X,\ad_Y^{2m}X]\in\Ls$. 

Now assume $n-m>1$ and write 
\begin{eqnarray*}
\lefteqn{\ad_Y[\ad_Y^{2n}X,\ad_Y^{2m}X]=}\\
 &&\qquad \ad_Y^2[\ad_Y^{2n-1}X,\ad_Y^{2m}X]
 -\ad_Y^2[\ad_Y^{2n-2}X,\ad_Y^{2m+1}X]+\ad_Y[\ad_Y^{2n-2}X,\ad_Y^{2m+2}X]. 
\end{eqnarray*}
The first two summands of the right hand side
belong to~$\Ls$ by the induction 
hypothesis on~$n+m$ 
and the last summand belongs to~$\Ls$ 
by the induction hypothesis on~$n-m$. 
Therefore~$\ad_Y[\ad_Y^{2n}X,\ad_Y^{2m}X]\in\Ls$. 

Next we show that 
$[\ad_Y^{2n}X,\ad_Y^{2m+1}X]\in\Ls$. For this purpose we
proceed by induction on~$n$. The case $n=0$ is precisely the 
assumption of the lemma. Now assume $n>0$ and write
\begin{eqnarray*}
\lefteqn{[\ad_Y^{2n}X,\ad_Y^{2m+1}X]=}\\
 &&\qquad \ad_Y^2[\ad_Y^{2n-2}X,\ad_Y^{2m+1}X]
 -2\ad_Y[\ad_Y^{2n-2}X,\ad_Y^{2m+2}X]+[\ad_Y^{2n-2}X,\ad_Y^{2m+3}X]. 
\end{eqnarray*}
The first summand of the right hand side belongs to $\Ls$ by the induction 
hypothesis on $n+m$, the second summand belongs to~$\Ls$ 
by the arguments in the previous paragraph,
and the last summand belongs to $\Ls$ 
by the induction hypothesis on $n$. Therefore 
$[\ad_Y^{2n}X,\ad_Y^{2m+1}X]\in\Ls$ and this completes 
the induction step on $n+m$. \EPf

\section{Applications}

In this section we describe some instances where 
the theorem can be applied. The most interesting case is probably
that of the complex hyperbolic plane~$\mathbf{H}^2(\mathbf C)$ 
which accomodates 
two classes of minimal hypersurfaces that are constructed by 
making use of the theorem. 

\subsection{Reflective submanifolds} An embedded submanifold $B$ of
a complete Riemannian manifold $M$ is called \emph{reflective}
if $B$ is complete with respect to the induced metric and it is a connected
component of the fixed point set of an involutive isometry of~$M$ 
(see~\cite{Le1}). By a well known result about fixed points of isometries, 
every reflective submanifold is automatically 
totally geodesic. 
In case $M$ is a simply connected globally symmetric space
and $\Lg=\Lk+\Lp$ is the associated involutive Lie algebra with 
respect to some fixed base point,
the reflective submanifolds $B$ through the base point
are in bijective correspondence with the subspaces~$\Lb$ of~$\Lp$ 
such that~$\Lb$ and its orthogonal complement~$\Lb^\perp$ 
are Lie triple systems and 
\[
 [[\Lb,\Lb^\perp],\Lb]\subset \Lb^\perp,\quad
[[\Lb,\Lb^\perp],\Lb^\perp]\subset\Lb. 
\]
Such a $\Lb$ is called a \emph{reflective subspace} of $\Lp$. 

If $B$ is a reflective submanifold of a nonpositively 
curved simply connected globally symmetric space $M$ and $v$ is 
\emph{any} nonzero normal vector to $B$,  
then it is readily seen that
the pair $(B,v)$ satisfies condition~(\ref{cond'}).
This shows that reflective submanifolds in nonpositively 
curved simply connected globally symmetric spaces 
may be translated along any normal geodesic to produce a 
minimal submanifold thereby providing an extensive number of 
situations to which our theorem can be applied. 
See~\cite{Le2,Le3} for tables with
the classification of reflective submanifolds in 
symmetric spaces. In the following we will comment 
on two interesting particular cases. 

\subsubsection{Rank one symmetric spaces} Every totally geodesic 
submanifold of a rank one globally symmetric space $M$ is reflective
(\cite{Le2,Le3}). 
We have that the maximal totally geodesic submanifolds 
of~$\mathbf{H}^n(\mathbf{C})$
are~$\mathbf{H}^{n-1}(\mathbf{C})$ and~$\mathbf{H}^n(\mathbf{R})$; 
of~$\mathbf{H}^n(\mathbf{H})$ (quaternionic hyperbolic space)
are~$\mathbf{H}^{n-1}(\mathbf{H})$ and~$\mathbf{H}^n(\mathbf{C})$; 
of~$\mathbf{H}^2(\mathbf{Cay})$ (Cayley hyperbolic plane)
are~$\mathbf{H}^2(\mathbf{H})$ and~$\mathbf{H}^8(\mathbf R)$ (\cite{Wo2}). 
In each case, the isotropy subgroup of the totally geodesic 
submanifold at the base point acts transitively on the unit sphere 
of the normal space at that point. 
Corresponding to each type above listed of totally geodesic 
submanifold in~$M$, 
the construction thus provides precisely one congruence
class of minimal embeddings~$\BR^N\to M$ irrespective of the 
choice of normal vector.

In particular, if~$M=\mathbf{H}^n(\mathbf{C})$ then the construction 
provides minimal embeddings~$\BR^{2n-1}\to\mathbf{H}^n(\mathbf{C})$ 
and~$\BR^{n+1}\to\mathbf{H}^n(\mathbf{C})$. 
The former submanifold is easily seen to be a bisector 
(see~\cite{Go}), whereas the latter one is a new example of minimal 
submanifold in~$\mathbf{H}^n(\mathbf{C})$ and it is a hypersurface if~$n=2$. 

\subsubsection{Real forms of Hermitian symmetric spaces}
A particularly interesting class of reflective submanifolds 
is composed of the fixed point sets of involutive antiholomorphic 
isometries of Hermitian symmetric spaces (which are automatically
connected, see~\cite{Le4}). Such reflective submanifolds
are called \emph{real forms} of the underlying 
Hermitian symmetric spaces, and correspond precisely to the
reflective subspaces $\Lb$ (cf.~\emph{supra})
which are totally real subspaces of $\Lp$.
The classification of real forms of Hermitian symmetric spaces
is given in~\cite{Le4}. It is interesting to note that 
we recover~$\mathbf{H}^n(\mathbf{R})$ as a real form 
of~$\mathbf{H}^n(\mathbf{C})$.

\subsection{Other examples} We now use the 
restricted root decomposition of a real semisimple Lie algebra
(see~\cite{He}, \S3, ch.~VI) 
in order to show that there are 
examples of totally geodesic submanifolds in symmetric spaces
which are not necessarily 
reflective but nevertheless can be translated along some 
normal geodesic to produce a minimal submanifold.
Let $\Lg=\Lk+\Lp$ be the involutive Lie algebra associated to a
symmetric space denote the involution by $\theta$. 
Choose a maximal Abelian subspace
$\La\subset\Lp$. Then there is an $\ad_{\La}$-invariant 
decomposition~$\Lg=\Lm+\La+\sum_{\lambda\in\Sigma}\Lg_\lambda$, where 
$\Lm$ is the centralizer of $\La$ in $\Lk$, $\Sigma$ 
is the restricted root system and, for each $\lambda\in\Sigma$,
$\Lg_\lambda$ is the corresponding restricted root space. 
Since $\theta[\Lg_\lambda]=\Lg_{-\lambda}$, we can write
\[ \Lk=\Lm+\sum_{\lambda\in\Sigma^+}\Lk_\lambda
\quad\mbox{and}\quad
\Lp=\La+\sum_{\lambda\in\Sigma^+}\Lp_\lambda, \]
where $\Lk_\lambda=\Lk\cap(\Lg_{\lambda}+\Lg_{-\lambda})$,
$\Lp_\lambda=\Lp\cap(\Lg_{\lambda}+\Lg_{-\lambda})$
and $\Sigma^+$ is the positive restricted root system. 
Now fix any $\lambda\in\Sigma^+$, define $\Ls=\Lp_\lambda$
and choose $X\in\La$. We claim that $\Ls$ is a Lie triple system 
and the pair $(\Ls,X)$ satisfies condition~(\ref{cond'}). 
This is to be a consequence of the following commutation rules 
(see~\cite{He}, Lemma~11.4, ch.~VII): 
\[ [\Lk_\lambda,\Lp_\mu]\subset\Lp_{\lambda+\mu}+\Lp_{\lambda-\mu},\quad
   [\Lk_\lambda,\Lk_\mu]\subset\Lk_{\lambda+\mu}+\Lk_{\lambda-\mu},\quad
   [\Lp_\lambda,\Lp_\mu]\subset\Lk_{\lambda+\mu}+\Lk_{\lambda-\mu}. \]
In fact, 
$[\Lp_\lambda,[\Lp_\lambda,\Lp_\lambda]]\subset
[\Lp_\lambda,\Lk_{2\lambda}+\Lm]\subset\Lp_\lambda$
so $\Lp_\lambda$ is a Lie triple system. 
Moreover, for any $Y\in\Lp_\lambda$ we have that $\ad_Y^{2n+1}X\in\Lk_\lambda$
and $\ad_Y^{2n}X\in\La+\Lp_{2\lambda}$ for $n=0,1,2,\ldots$ 
Therefore $[X,\ad_Y^{2n}X]\in\Lp_{\lambda}$ for $n=0,1,2,\ldots$

\bibliographystyle{amsplain}
\bibliography{ref}

\end{document}